\documentclass[11pt]{article}

\usepackage[bookmarks=false]{hyperref}
\hypersetup{colorlinks,
      linkcolor=blue,
      citecolor=blue,
      urlcolor=blue}
\usepackage[figuresright]{rotating}
\usepackage[utf8]{inputenc} 
\usepackage[T1]{fontenc}    
\usepackage{amsmath}
\usepackage{epsfig,xcolor}
\usepackage[tworuled,linesnumbered,vlined]{algorithm2e} 
\usepackage{wasysym}
\usepackage{esvect}
\usepackage{amsthm}
\usepackage{amssymb}

\newtheorem{theorem}{Theorem}[section]
\newtheorem{definition}[theorem]{Definition}
\newtheorem{lemma}[theorem]{Lemma}
\newtheorem{proposition}[theorem]{Proposition}

\newcommand{\dd}{\ensuremath{\vv{D}}}
\newcommand{\dt}{\vv{T}\!}
\newcommand{\A}{\ensuremath{\mathcal{A}}\xspace}

\newcommand{\Prime}{\mathsf{P}}
\newcommand{\Komplete}{\mathsf{K}}
\newcommand{\Star}{\mathsf{S}}
\newcommand{\X}{\mathsf{X}}

\newcommand{\algo}{\mathcal{A}\xspace}

\newcommand{\gp}{\mathit{gp}}

\definecolor{linecomment}{rgb}{0.95, 0.1, 0.1}

\begin{document}

\title{Mutual-visibility in distance-hereditary graphs: \\ a linear-time algorithm}
\date{}

\author{
Serafino Cicerone $^{a}$\thanks{Email: \texttt{serafino.cicerone@univaq.it}}
\and
Gabriele Di Stefano $^{a}$\thanks{Email: \texttt{gabriele.distefano@univaq.it}}
}

\maketitle
\begin{center}
\small
$^a$ Department of Information Engineering, Computer Science, 
and Mathematics, \\	     
University of L'Aquila, Italy 
\end{center}

\begin{abstract}
The concept of mutual-visibility in graphs has been recently introduced. If $X$ is a subset of vertices of a graph $G$, then vertices $u$ and $v$ are $X$-visible if there exists a shortest $u,v$-path $P$ such that $V(P)\cap X \subseteq \{u, v\}$. If every two vertices from $X$ are $X$-visible, then $X$ is a mutual-visibility set. The mutual-visibility number of $G$ is the cardinality of a largest mutual-visibility set of $G$. It is known that computing the mutual-visibility number of a graph is NP-complete, whereas it has been shown that there are exact formulas for special graph classes like paths, cycles, blocks, cographs, and grids. In this paper, we study the mutual-visibility in distance-hereditary graphs and show that the mutual-visibility number can be computed in linear time for this class. 
\end{abstract}


\maketitle

%
\section{Introduction}
Given a set of points in Euclidean space, they are mutually visible if and only if three of them are not collinear. In other words, two points $p$ and $q$ are mutually visible when no other points belong to the segment $pq$. A line segment in Euclidean space represents the shortest path between two points, but in more general topologies, this type of path (called geodesic) may not be unique. Then, in general, two points are mutually visible when there exists at least a shortest path between them without further points. 

In~\cite{DiStefano22}, this concept has been recently extended to mutual-visibility in graphs. Formally, given a graph $G$ and a set of vertices $X\subseteq V(G)$, two vertices $u$ and $v$ are $X$-visible if there exists a shortest $u,v$-path $P$ such that $V(P)\cap X \subseteq \{u,v\}$. If every two vertices from $X$ are $X$-visible, then $X$ is a mutual-visibility set of $G$. The mutual-visibility number of $G$, denoted as $\mu(G)$, is the cardinality of a largest mutual-visibility set of $G$. In the same original paper, it is shown that computing $\mu(G)$ is a NP-complete problem. Moreover, it has been shown that there are exact formulas for special graph classes like paths, cycles, blocks, cographs, and grids. In~\cite{CiceroneDK23,CiceroneDKY23}, exact formulas have been derived also for both the Cartesian and the Strong product of graphs.

In this paper, we study the mutual-visibility in distance-hereditary graphs. These graphs have been introduced by Howorka in~\cite{howorka:77}, and are defined as those graphs in which every connected induced subgraph is isometric, that is the distance between any two vertices in the subgraph is equal to the one in the whole graph. Therefore, any connected induced subgraph of any distance-hereditary graph $G$ ``inherits'' its distance function from $G$. This kind of graphs have been rediscovered many times (e.g., see~\cite{bandelt/mulder:86,peerj21,dam01}). Since their introduction, dozens of papers have been devoted to them, and different kind of characterizations have been found: metric, forbidden subgraphs, cycle/chord conditions, level/neighborhood conditions, generative, and more (e.g., see~\cite{BrandstadtLS99}). Note that distance-hereditary graphs include some classes for which the mutual-visibility problem has been already studied: trees, block graphs, and cographs. 

As the main result, in this paper, we show that the mutual-visibility number can be computed in linear time for distance-hereditary graphs. This result exploits the well-known characterization based on the split decomposition~\cite{bouchet:88}.

\smallskip\noindent
\textbf{Motivations and related works.} 
%
 Vertices in mutual visibility may represent entities on some nodes of a computer/social network that want to communicate in a efficient and ``confidential'' way, that is, in such a way that the exchanged messages  do not pass through other entities. Moreover, in the context of swarm robotics, mobile robots or agents usually move and operate in a continuous environment (e.g., the Euclidean plane) or in a discrete space (e.g., a graph)~\cite{FPS19}. In many cases,  it is required to locate (or move) the agents in positions so that they are all in mutual visibility (e.g., see~\cite{LunaFCPSV17,PoudelAS21}). Note that in these works, even if the agents are constrained to move in a grid, the visibility is always checked along line segments. A first result in which agents operate in graphs and the visibility is checked along shortest paths can be found in~\cite{CiceroneFSN23}.

A concept similar to the mutual-visibility is given by the \emph{general position}. It has been introduced in~\cite{ManuelK18} and defined as follows: in a graph $G$, a general position set is a subset $S\subseteq V(G)$ such that no three vertices from $S$ lie in a common geodesic of $G$. The \emph{general position number} $\gp(G)$ is the order of a largest general position set of $G$.
Since its introduction, the general position number has been studied for several graph classes (e.g., grid networks, cographs and bipartite graphs, Cartesian products of graphs, and Kneser graphs.
The difference between a general position set $S$ and a mutual-visibility set $X$ is that two vertices are in $X$ if there \emph{exists} a shortest path between them with no further vertex in $X$, whereas two vertices are in $S$ if for \emph{every} shortest path between them, no further vertex is in $S$. The two concepts are intrinsically different, but also closely related, since the vertices of a general position set are in mutual visibility.  


%
\section{Notation and preliminaries}\label{sec:notation}
In this work, 
given a graph $G$, $V(G)$ and $E(G)$ are used to denote its vertex set and its edge set, respectively. We use standard terminologies from~\cite{BrandstadtLS99}, some of which are briefly reviewed here.

If $X\subseteq V(G)$, then $G[X]$ denotes the subgraph of $G$ induced by $X$, that is the maximal subgraph of $G$ with vertex set $X$. $N_G(X)$ is the open neighborhood of $X$ in $G$, that is $N_G(X) = \{v\in V(G)\setminus X:~ \exists u\in X, uv\in E(G)\}$. If $X = \{u\}$ we simply write $N_G(u)$ instead of $N_G(\{u\})$. If $|N_G(u)|=1$, $u$ is called \emph{pendant} vertex.

%

We call $G_1$ and $G_2$ \emph{isomorphic}, and write $G_1\sim G_2$, if there exists a bijection $\varphi :V(G_1) \rightarrow V(G_2)$ such that $uv \in E(G_1) \Leftrightarrow \varphi(u)\varphi(v) \in E(G_2)$ for all $u,v \in V(G_1)$. Such a bijection $\varphi$ is called \emph{isomorphism}.

An \emph{edge cut-set} of a connected graph $G$ is any subset $E'\subseteq E(G)$ of edges that, if removed, forms a graph with more than one connected component. An edge cut-set of size 1 is simply called \emph{cut-edge}.
Similarly, a \emph{vertex cut-set} is any subset $V'\subseteq V(G)$ such that the induced subgraph $G[V(G)\setminus V']$ has more than one connected component. A vertex cut-set of size 1 is called \emph{cut-vertex}.
A graph $G$ is \emph{biconnected} if it has no cut-vertices. $G$ is a \emph{block graph} if each biconnected component (``block'') of $G$ is a clique.
The \emph{complete graph} (or \emph{clique}) $K_n$, $n\ge 1$, is the graph with $n$ vertices where each pair of distinct vertices are adjacent. The \emph{path graph} $P_n$, $n\ge 2$, is the graph with $V(P_n) = \{v_1,v_2,\ldots,v_n\}$ such that $v_i$ is adjacent to $v_j$ if and only if $|i-j|=1$. A \emph{complete bipartite graph} $K_{m,n}$ is a graph whose vertices can be partitioned into two subsets $V_1$ and $V_2$, with $|V_1|=m$ and $|V_2|=n$, such that no edge has both endpoints in the same subset, and every possible edge that could connect vertices in different subsets is part of the graph. A \emph{star graph} corresponds to any $K_{1,n}$.
%

%


Concerning the mutual-visibility, it is easy to observe that $\mu(G)\ge 1$ for each graph $G$ (indeed, any vertex $u\in V(G)$ is a mutual-visibility set of $G$). 

From~\cite{DiStefano22} we also know that computing the mutual-visibility number of a graph is NP-complete. Concerning small values of $\mu$, we know that:
\begin{itemize}
\item $\mu(G)=1$ if and only if $G \sim K_1$; 
\item $\mu(G)=2$ if and only if $G \sim P_n$, $n\ge 2$;
\end{itemize}
Moreover, a partial characterization for $\mu(G)=3$ is provided in~\cite{CiceroneDK23}. 

The following lemma asserts that there always exists a maximum mutual-visibility set for a graph $G$ without cut-vertices.
\begin{lemma} \emph{ (\cite[Lemma 2.5]{DiStefano22}) } 
\label{lem:cut} 
Let $C$ be the set of the cut-vertices in a graph $G$. There exists a $\mu$-set $X$ for $G$ such that $X \cap C = \emptyset$.
\end{lemma}



\begin{figure}[h]
\begin{center}
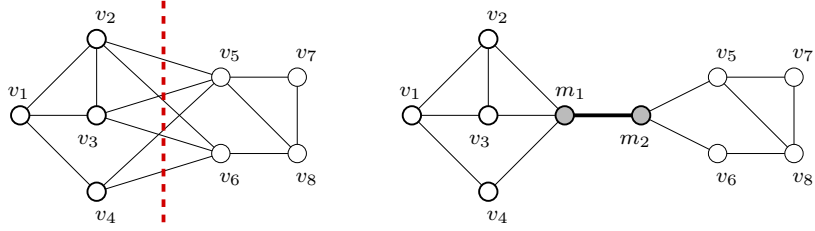
\end{center}
\caption{\em A graph $G$ along with the marked graph forming a simple decomposition of $G$. Grey vertices and bold lines are used to visualize marked vertices and marked edge, respectively. Notice that each obtained bag is isomorphic to an induced subgraph of $G$.}
\label{fig:split}
\end{figure}

We recall the definition of split decomposition given in~\cite{bouchet:88}. A split in a graph $G$ is a vertex partition $(X,Y)$ of $G$ such that $|X|,|Y| \ge 2$ and every vertex of $N_G(X)$ is adjacent to every vertex of $N_G(Y)$. Notice that not all connected graphs have a split, and those that do not have a split are called \emph{prime graphs}.
A \emph{marked graph} $D$ is a connected graph with a set of edges $M(D)$, called 
\emph{marked edges}, that form a matching such that every edge in $M(D)$ is a cut-edge. 
The ends of the marked edges are called \emph{marked vertices}, and the components of $
(V(D),E(D)\setminus M(D))$ are called \emph{bags} of $D$. The edges in $E(D)\setminus 
M(D)$ are called \emph{unmarked edges}, and the vertices that are not marked vertices are 
called \emph{unmarked vertices}. If $(X,Y)$ is a split in $G$, then we construct a marked 
graph $D$ that consists of the vertex set $V(G) \cup \{x', y'\}$ for two distinct new 
vertices $x', y' \not\in  V(G)$ and the edge set 
$E(G[X]) \cup E(G[Y]) \cup \{x'y'\} \cup E'$, where we define $x'y'$ as marked and
$ E' := \{x'x ~|~ x \in X \wedge \exists  y\in Y: xy \in E(G) \} \cup 
\{y'y ~|~ y \in Y  \wedge \exists x\in X: xy \in E(G) \}.$
The marked graph $D$ is called a \emph{simple decomposition} of $G$ (e.g., see Figure~\ref{fig:split}). 

The next proposition follows from the definition of split.
\begin{proposition}\label{prop:cut}
Let $(X,Y)$ be a split of a graph $G$.  If $X\setminus N_G(Y)$ ($Y\setminus N_G(X)$, resp.) is not empty, then $N_G(Y)$ ($N_G(X)$, resp.) is a vertex cut-set of $G$.
\end{proposition}
%
%

A \emph{split decomposition} of a graph G is a marked graph $D$ defined inductively to be either $G$ or a marked graph defined from a split decomposition $D'$ of $G$ by replacing a connected component $H$ of $(V (D'), E(D')\setminus M(D'))$ with a simple decomposition of $H$. 

For a marked edge $xy$ connecting two bags $X$ and $Y$ in a split decomposition $D$, the \emph{recomposition} of $D$ along $xy$ is the split decomposition $D'$ obtained by making each vertex in $N_{X}(x')$ adjacent to each vertex in $N_{Y}(y')$ and by removing $x'$, $y'$, and their adjacent edges. For a split decomposition $D$, let $G[D]$ denote the graph obtained from $D$ by recomposing all marked edges. By definition, if $D$ is a split decomposition of $G$, then $G[D] = G$. Since each marked edge of a split decomposition $D$ is a cut-edge and all marked edges form a matching, if we contract all unmarked edges in $D$, then we obtain a tree. We call it the \emph{decomposition tree} of $G$ associated with $D$ and denote it by $T(D)$. To distinguish the vertices of $T(D)$ from the vertices of $G$ or $D$, the vertices of $T(D)$ will be called \emph{nodes}. Obviously, the nodes of $T(D)$ are in bijection with the bags of $D$. Two bags of $D$ are called \emph{neighbour bags} if their corresponding nodes in $T(D)$ are adjacent. 

\begin{figure}[t]
\begin{center}
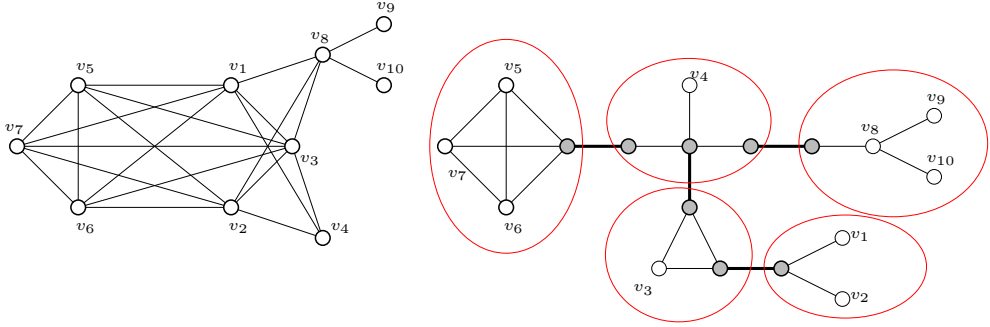 
\end{center}
\caption{\em A distance-hereditary graph $G$ along with $D_G$. 
Red ovals, grey vertices, and bold lines are used to visualize the bags, the marked vertices, and the marked edges of $D_G$.}
\label{fig:G}
\end{figure}

A split decomposition $D$ of $G$ is called a \emph{canonical split decomposition} (or canonical decomposition for short) if each bag of $D$ is either a prime graph, a star, or a complete graph, and $D$ is not the refinement of a decomposition with the same property. Figure~\ref{fig:G} shows a graph $G$ along with its canonical decomposition.

\begin{theorem}\emph{(Cunningham and Edmonds~\cite{cunningham/edmonds:80}, Dahlhaus~\cite{dahlhaus:jagm00})}
\label{teo:dahlhaus_complexity}
Every connected graph $G$ has a unique canonical decomposition, up to isomorphism, and it can be computed in time $O(|V(G)| + |E(G)|)$.
\end{theorem}

From Theorem~\ref{teo:dahlhaus_complexity}, we can talk about only one canonical decomposition of a graph G because all canonical decompositions of $G$ are isomorphic. Let $D$ be a split decomposition of a graph $G$ with bags that are either prime graphs, complete graphs or stars. 
The type of a bag of $D$ is either $\Prime$, $\Komplete$, or $\Star$ depending on whether it is a prime graph, a complete graph, or a star. The type of a marked edge $uv$ is $AB$ where $A$ and $B$ are the types of bags containing $u$ and $v$ respectively. If $A = \Star$ or $B = \Star$, then we can replace $\Star$ with $\Star_p$ or $\Star_c$ depending on whether the end of the marked edge is a pendant vertex (i.e., a leaf) or the center of the star.
%
%
%
%
%
We denote by $D_G$ the canonical decomposition of $G$, by $M_G$ the marked edges of $D_G$, and by $T_G$ the decomposition tree associated with $D_G$. 
We remark on the well-known property for which each bag of $D_G$ is isomorphic to an induced subgraph of $G$. 
We use the next characterization of distance-hereditary graphs. 

\begin{theorem}\emph{(Bouchet~\cite{bouchet:88})}
\label{teo:teo:bouchet_canonical_DH}
A graph $G$ is distance-hereditary if and only if each bag of $D_G$ is of type $\Komplete$ or $\Star$ and no marked edge of $D_G$ is of type $\Komplete\Komplete$ or $\Star_p\Star_c$.
\end{theorem}

According to this theorem, if $G$ is distance-hereditary we denote any bag of $D_G$ as $\Komplete$-bag or $\Star$-bag. An \emph{alternating path} is any path connecting two  
vertices in $D_G$ such that an unmarked edge and a marked edge alternatively appear in the path (e.g., the path from $v_1$ to $v_7$ in the canonical decomposition shown in Figure~\ref{fig:G}). 

\begin{lemma}\label{lem:alternating} \emph{(Courcelle~\cite{Courcelle06})}
Let $G$ be a graph. Then, $xy\in E(G)$ if and only if there exists an alternating path between $x$ and $y$ in $D_G$. This alternating path is moreover unique.
\end{lemma}

We conclude this section by providing a useful statement. 

\begin{lemma}\label{lem:star}
Let $G$ be a graph. 
If $v\in V(G)$ is the center of a $\Star$-bag of $D_G$ then it is a cut-vertex of $G$. 
\end{lemma}
\begin{proof}
If $D_G$ has only one bag, according to the hypothesis $G$ is a star with $v$ as the central vertex. In this case, the statement trivially follows. Assume now that $D_G$ has at least two bags and let $X$ be the set of vertices of $D_G$ forming the $\Star$-bag of which $v$ is the center. Let $u_1,u_2\in X\setminus \{v\}$. Since $D_G$ has at least two bags, we can assume that $u_1$ is marked. Consider any alternating path $v,u_1,w_1,w_2, \ldots, w_k$, with $k\ge 2$ and $w_k\in V(G)$. Since $u_1w_1$ is a marked edge, by Proposition~\ref{prop:cut} we easily get that removing $v$ in $G$ would disconnect $w_k$ from either $u_2$ (if $u_2$ is unmarked) or $u_2'$ (if $u_2$ is marked and $u_2,\ldots,u_2'$ is any alternating path from $u_2$ to an unmarked vertex $u_2'$).
\end{proof}

According to this property, we call \emph{$\sigma$-vertex} any vertex $v\in V(G)$ which is the center of a $\Star$-bag. 
%

%
\section{The directed canonical decomposition}\label{sec:main}

%
%
We introduce the notion of \emph{directed canonical decomposition} of a distance hereditary graph $G$ useful for the construction of the algorithm capable of computing a $\mu$-set of $G$, described in the next section.

\begin{definition}\label{def:dd}
Let $G$ be a distance hereditary graph. 
\begin{itemize}
    \item 
    The \emph{directed canonical decomposition} of $G$ is denoted as $\dd_G$, and corresponds to the mixed graph obtained from $D_G$ according to the following edge-replacing operation: for each marked edge $uv$ of type $\Star_c\X$, with $\X\in \{\Star_c,\Komplete\}$, remove $uv$ from $M_G$ and insert an edge directed from $u$ to $v$ into $M_G$.
    \item
    A directed decomposition tree $\dt_G$ can be associated with $\dd_G$ in the same way $T_G$ has been associated with the canonical decomposition $D_G$.
\end{itemize}  
\end{definition}


\begin{figure}[t]
\begin{center}
\includegraphics[width=12.0cm]{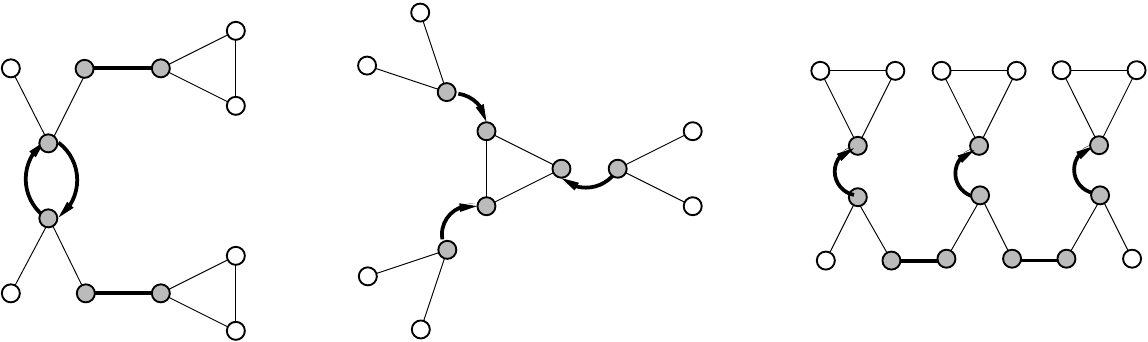}
\caption{\em Three examples of directed canonical decomposition. In order, two opposite t-arrows, head-connected t-arrows, and tail-connected t-arrows.
}
\label{fig:cases}
\end{center}
\end{figure}

According to Definition~\ref{def:dd}, when forming $\dd_G$, a marked edge $uv$ can be replaced by one oriented edge (when $uv$ is of type $\Star_c\Komplete$) or by two opposite oriented edges (when $uv$ is of type $\Star_c\Star_c$). For sake of simplicity, we call \emph{arrow} an oriented marked edge $a=(u,v)$, and we call \emph{opposite} two distinct arrows joining the same pair of vertices. 
Given an arrow $a=(u,v)$, we call $u$ ($v$, resp.) \emph{tail} (\emph{head}, resp.) of $a$.

Given an arrow $a=(u,v)$, we denote by $\dd_G^t(a)$ and $\dd^{h}_G(a)$ the two subgraphs of $\dd_G$ obtained by removing $a$ (and the arrow $(v,u)$ if present) containing the tail $u$ and the head $v$, respectively. 
Similarly, we denote as $V_G^h(a)$ ($V_G^t(a)$, resp.) the set containing all the unmarked vertices in $\dd_G^h(a)$ ($\dd_G^t(a))$, resp.) reachable from $v$ ($u$, resp.) via an alternating path.  

\begin{definition}
    Given a distance-hereditary graph $G$, an arrow $a=(u,v)$ in $\dd_G$ is called \emph{terminal} (t-arrow, for short) if both the following conditions hold:
    \begin{enumerate}
        \item $V_G^h(a)$ does not contain a $\sigma$-vertex;
        \item $\dd_G^h(a)$ does not contain an alternating 
                 $v,u'$-path in $D_G$ followed by an arrow $(u',v')$.
    \end{enumerate}
\end{definition}

\begin{figure}[t]
\begin{center}
\includegraphics[width=11.5cm]{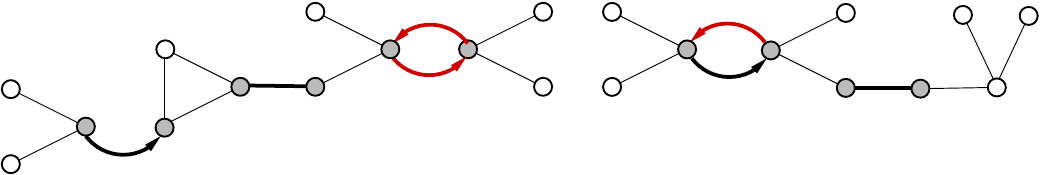}
\caption{\em Two directed canonical decompositions (t-arrows are represented in red).
}
\label{fig:t-arrow}
\end{center}
\end{figure}


%
The following statements provide some structural properties of t-arrows.%

\begin{lemma}\label{lem:sym-arrows}
Given a distance-hereditary graph $G$, if $a=(x,y)$ and $b=(y,x)$ are opposite t-arrows in $\dd_G$, then they are the only t-arrows in $\dd_G$.
\end{lemma}
\begin{proof}
Assume that there is an arrow $e=(u,v)$ in $\dd_G^h(a)=\dd_G^{t}(b)$.
On one hand, since $a$ is a t-arrow, there is no alternating path connecting the marked vertex $y$ to the marked vertex $u$, tail of $e$. Since $b$ is a t-arrow, there is no alternating path connecting the marked vertex $y$ to the marked vertex $v$, head of $e$.
On the other hand, any shortest path $P$ connecting $y$ to $u$ or to $v$, that is not alternating,  must contain two consecutive unmarked edges (two consecutive marked edges are not possible). 
Consider the first two consecutive unmarked edges encountered on $P$ starting from $y$. These two edges cannot belong to a $\Komplete$-bag otherwise $P$ would not be a shortest path. Hence, the two edges must belong to a $\Star$-bag. If the center of this $\Star$-bag is a marked vertex, then the corresponding marked edge is an arrow, but it is not possible since $a$ is a t-arrow. If the center of the $\Star$-bag is unmarked, then it is a $\sigma$-vertex and $a$ would not be a t-arrow, by definition. So $\dd_G^h(a)$ does not contain any t-arrow. 
By symmetry, also $\dd_G^h(b)=\dd_G^{t}(a)$ does not contain t-arrows.
\end{proof}

The following three lemmas show how the t-arrows appear in 
$\dd_G$.

\begin{lemma}\label{lem:out-arrows}
If $G$ is a distance-hereditary graph, then in $\dd_G$ there is no path $x,y,u_1,u_2,\ldots,u_k, x',y'$ in which both $(x,y)$ and $(x',y')$ are t-arrows. 
\end{lemma}
\begin{proof}
Assume there exists in $\dd_G$ a path $x,y,u_1,u_2,\ldots,u_k,x',y'$ in which both $(x,y)$ and $(x',y')$ are t-arrows. By the definition of t-arrow, it cannot be an alternating path. 
Then, consider any shortest path $P$ defined as $x,y,u'_1, u'_2, \ldots, u'_{k'}, x',y'$. Note that $k'\ge 1$ otherwise $P$ reduces to $x,y,x',y'$, an alternating path. As shown in the proof of Lemma~\ref{lem:sym-arrows}, $P$ contains two consecutive unmarked edges belonging to a $\Star$-bag.  In the same proof, it is shown that the presence of this $\Star$-bag implies that $(x,y)$ cannot be t-arrow, a contradiction.
\end{proof}

Assuming that $\dd_G$ has more than one t-arrow, then either Lemma~\ref{lem:sym-arrows} holds and there are only two opposite t-arrows, or, by Lemma~\ref{lem:out-arrows}, a shortest path connecting the vertices of two t-arrows $(x,y)$ and $(x',y')$ is either in the form $x,y=u_0,u_1,u_2,\ldots,u_k=y',x'$ or in the form $y,x=u_0,u_1,u_2,\ldots,u_k=x',y'$, $k\geq 1$. If the path is in the first form, we say that the arrows are \emph{head-connected}, otherwise \emph{tail-connected}. 

 
\begin{lemma}\label{lem:head-tail}
Let $G$ be a distance-hereditary graph. If $\dd_G$ has more than one t-arrow, they are 
two and opposite, or pairwise head-connected, or pairwise tail-connected.
\end{lemma}
\begin{proof}
The first case follows directly from Lemma~\ref{lem:sym-arrows}. Suppose now that in $\dd_G$ there are not two opposite t-arrows. By contradiction assume that all the t-arrows in $\dd_G$ are neither pairwise head-connected nor pairwise tail-connected. Then there must exist three t-arrows $a=(x,y)$,  $b=(x',y')$, and $c=(x'',y'')$ such that two of them, say $a$ and  $b$, are head-connected and the third t-arrow $c$ that is tail-connected with $a$ or $b$. Without loss of generality, assume that $c$ is tail-connected with $a$. Then, there must exist a path $y,x,u_1,\ldots,u_k,x'',y''$. Hence, the path $x',y',v_1,\ldots,v_{k'},y,x,u_1,\ldots,u_k,x'',y''$ contradicts Lemma~\ref{lem:out-arrows}. 
\end{proof}


Interestingly, For a distance-hereditary $G$, if the t-arrows of $\dd_G$ are pairwise head-connected they are in a particular configuration.

\begin{lemma}\label{lem:clique}
Let $G$ be a distance-hereditary graph such that all the pairs of t-arrows in $\dd_G$ are head-connected.  Then, the heads of all the t-arrows belong to a unique $\Komplete$-bag.
\end{lemma}
\begin{proof}
Let $a=(x,y)$ and $b=(x',y')$ be two t-arrows. Assume that the shortest path $P$ connecting $y$ and $y'$ contains an edge $e$ from a $\Star$-bag. If $P$ is an alternating path, one of the two vertices of $e$ is the marked center of the $\Star$-bag. The corresponding marked edge is an arrow: impossible because one between $a$ and $b$ would not be a t-arrow. If $P$ is not an alternating path there are two edges, say $e$ and $e'$, of a $\Star$-bag in $P$ that are both incident to the center of the $\Star$-bag. Whether or not the center is a $\sigma$-vertex, both $a$ and $b$ would not be t-arrows. Then, $P$ is an alternating path and contains only edges from $\Komplete$-bags. By Theorem~\ref{teo:teo:bouchet_canonical_DH},  there cannot be two adjacent $\Komplete$-bags and then the $\Komplete$-bag is unique and $P$ consists only of one edge from this clique. 
\end{proof}

\begin{figure}[t]
\begin{center}
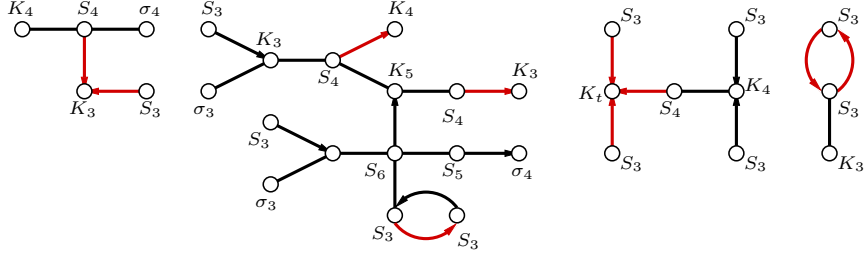
\end{center}
\caption{\em Four directed decomposition trees named, in order, $\protect\vv{T}\!_G$,
$\protect\vv{T}\!_{G_1}$, 
$\protect\vv{T}\!_{G_2}$, and $\protect\vv{T}\!_{G_3}$. Each vertex is labelled with $K_n$, $S_n$, or $\sigma_n$ to recall the corresponding bag in $\protect\vv{D}_{G}$ (a clique, a star without $\sigma$ vertex, a start with $\sigma$ vertex). The subscript $n$ denotes the size of that bag. For instance, $\protect\vv{T}\!_G$ corresponds to the directed canonical decomposition shown in Figure~\ref{fig:G}.}
\label{fig:dts}
\end{figure}

\begin{lemma}\label{lem:block}
     Given a distance-hereditary graph $G$, $\dd_G$ has no arrows if and only if $G$ is a block graph. 
\end{lemma}
\begin{proof}
($\Leftarrow$) Assuming that $a$ is an arrow of $\dd_G$, recomposing the marked edge in $D_G$ corresponding to $a$ would lead to a biconnected subgraph of $G$ that is not a clique, contradicting the definition of block graph.

($\Rightarrow$) In this case all the marked edges in $D_G$ are of type $\Komplete\Star_p$ or $\Star_p\Star_p$ whose recompositions lead to biconnected subgraphs of $G$ that are cliques: $K_2$ in the latter case, and $K_n$, $n\ge 3$, in the former.
\end{proof}

%
The following statements motivate the introduction of t-arrows with respect to the need to compute a $\mu$-set of a distance-hereditary graph $G$.

\begin{lemma}\label{lem:arrow}
    Given a distance-hereditary graph $G$, a set $X\subseteq V(G)$ and an arrow $a$ in  $\dd_G$, two vertices $x,y\in V_G^t(a)$ are $X$-visible if there exist $w\in V_G^h(a)\setminus X$.
\end{lemma}
\begin{proof}
   Since there are two alternating paths from $x$ to $w$ and from $w$ to $y$ in $\dd_G$,  then $w\in N_G(x)\cap N_G(y)$. This implies that $x$ and $y$ are at distance at most two. If $x$ and $y$ are adjacent, then are obviously $X$-visible. If they are at distance two, path $x,w,y$ is a shortest path, and since $w$ is not in $X$ vertices $x$ and $y$ are $X$-visible. 
\end{proof}


\begin{lemma}\label{lem:2arrows}
     Given a distance-hereditary graph $G$, a set $X\subseteq V(G)$ and two arrows $a=(u,v)$ and $b=(u',v')$ in  $\dd_G$ connected by an alternating path $u,v=v_0,v_1,\ldots, v_k=u',v'$, two vertices $x,y\in V_G^t(a)$ are $X$-visible if there exists $w\in V_G^h(b)\setminus X$.
\end{lemma}
\begin{proof}
    The proof reduces to Lemma~\ref{lem:arrow} by simply observing that $V_G^t(a)\subseteq V_G^t(b)$ and $V_G^h(a)\supseteq V_G^h(b)$.
\end{proof}

\begin{algorithm}[h!]
\smaller
	\SetKwInput{Proc}{Algorithm \A}
	\Proc{}
	\SetKwInOut{Input}{Input}
	\Input{a distance-hereditary graph $G$.}
	\SetKwInOut{Output}{Output}
	\Output{a $\mu$-set $X$ of G.} 
	
	\BlankLine
	\BlankLine
    Compute $\dd_G$\;
    Let $S$ be the set of the $\sigma$-vertices of $\dd_G$\;
    Let $X =  V(G)\setminus S$\;
    Let $A$ be a set containing all the t-arrows of $\dd_G$\;
    \BlankLine
    
    \lIf{$A=\emptyset$}{\Return $X$}\label{line:Avuoto}


    \If{$A$ contains either a single t-arrow or pairwise tail-connected t-arrows  \label{line:tails}} {
        \ForAll{ t-arrow $a$ in $A$} {
            Let $w\in V_G^h(a)$\;
            $X := X\setminus \{w\}$
        }
        \Return $X$  \label{line:tails-end}
    }

    \If{$A$ contains pairwise head-connected t-arrows\label{line:heads}} {
        Let $K$ be the $\Komplete$-bag containing all the heads of the t-arrows in $A$\;
        \lIf{$K$ has an unmarked vertex $w$\label{line:K}} {
          $X := X\setminus \{w\}$
        }
        \Else{
        Let $C:=\{\dd^t_G(a)~|~ a\in A\}$\;\label{line:C} 
        \lIf{$C$ contains a special element $D$}{
          $X := X\setminus \{w\}$, with $w$ special vertex of  $D$ \label{line:one-special} 
        }
        \Else{
          Let $D$ and $D'$ be two elements of $C$\;
            Let $w\in V(D)$ and $w'\in V(D')$ be unmarked vertices\; \label{line:two-vertices}
            $X := X\setminus \{w,w'\}$  \label{line:end-special1}\;
        }
        }
        \Return $X$ \label{line:heads-end}
    }
      
    \If{$A$ contains two opposite t-arrows $a$ and $b$ \label{line:opposite}}{ 
    
    \lIf{$\dd^h_G(a)$ is special}
                {$X := X\setminus \{w\}$, with $w \in V_G^h(a)$ special vertex of $\dd^h_G(a)$\label{line:remove1}}
    \Else{
         \lIf{$\dd^h_G(b)$ is special }
                {$X := X\setminus \{w'\}$, with $w' \in V_G^h(b)$ special vertex of $\dd^h_G(b)$\label{line:remove2}}
         \Else{
            Let $w\in V_{G}^t(a)$ and $w'\in V_{G}^t(b)$  be unmarked vertices\; \label{line:two-vertices2} 
            $X := X\setminus \{w,w'\}$\label{line:end-special2}
         }   
    }
    \Return $X$  \label{line:opposite-end}
    }
\caption{Algorithm for computing a $\mu$-set $X$ of any distance-hereditary graph $G$. }
	\label{alg:A}
\end{algorithm}

\section{The algorithm}\label{sec:algorithm}
%
%
We have already remarked that the directed canonical decomposition of a distance hereditary graph G has been formulated to be a valid tool for computing a $\mu$-set of $G$. The algorithm we propose for this purpose is simply denoted as $\algo$ and its pseudocode can be fond in Algorithm~\ref{alg:A}. Basically, it first computes $\dd_G$ and then performs specific actions according to the number and the structure of the t-arrows of $\dd_G$. In particular, the strategy underlying $\algo$ is based on the following steps, in order: 
\begin{itemize}
    \item it is initially assumed that the entire $V(G)$ can play the role of a $\mu$-set denoted as $X$;
    \item then, if there are $\sigma$ vertices in $\dd_G$, they are removed from $X$ (indeed, Lemma~\ref{lem:star} states that $\sigma$ vertices are cut-vertices in $G$, and hence, by Lemma~\ref{lem:cut}, they can be correctly eliminated from $X$);
    \item finally, t-arrows of $\dd_G$ are analyzed: 
    \begin{itemize}
        \item  if there are no arrows at all (in this case, Lemma~\ref{lem:block} states that $G$ is a block graph), or no t-arrows, then it is possible to prove that the computed set $X$ is indeed a $\mu$-set of $G$;
        \item  if there are t-arrows, by Lemmas~\ref{lem:arrow} and~\ref{lem:2arrows} some vertices must be further removed from $X$. 
        For such cases, $\algo$ follows three different approaches according to the structural characterization given by Lemma~\ref{lem:head-tail}.
    \end{itemize}
\end{itemize}

To describe the algorithm in terms of these three approaches, in what follows we use the examples illustrated by the directed decomposition trees named $\dt_{G_1}$, $\dt_{G_2}$, and $\dt_{G_3}$ and shown in Figure~\ref{fig:dts}.

Tree $\dt_{G_1}$ is obtained from a graph $G_1$ in which its directed canonical decomposition has many t-arrows that are pairwise tail-connected. In this case, $\algo$ first removes from $X=V(G)$ three $\sigma$ vertices and then executes the block of Lines~\ref{line:tails}--\ref{line:tails-end}. In such a block, for each t-arrow $a$, an arbitrary vertex $w\in V_{G_1}^h(a)$ is removed from the current set $X$. After all the t-arrows are processed, $X$ is returned as a $\mu$-set of $G$.

In the second example, $\dt_{G_2}$ is obtained from a graph $G_2$ in which its directed canonical decomposition has pairwise head-connected t-arrows. In this case, $\algo$ executes the block of Lines~\ref{line:heads}--\ref{line:heads-end}. According to Lemma~\ref{lem:clique}, the heads of all the t-arrows belong to a unique $\Komplete$-bag (see $K_t$ in the example). If $K_t$ has an unmarked vertex $w$ (in the example, $K_t$ should have $t\ge 4$ vertices), then $X\setminus \{w\}$ is returned. This is because, for each arrow $a$ in $\dd_{G_2}$, each pair of elements $u,v\in V_{G_2}^t(a)$ are connected to $w$, and the removal of $w$ from $X$ makes them $X$-visible. However, if $w$ does not exist, that is when all the vertices of $K_t$ are marked, still there are cases in which the removal of a single vertex is sufficient. These cases occur when a generic directed canonical decomposition $\dd_G$ has a ``special subgraph'' as in the following definition. 

\begin{definition}
   We call \emph{special} any subgraph $D$ of $\dd_G$ isomorphic to a single $\Star$-bag with a marked center and (1) two unmarked pendant vertices $x$ and $y$,  or (2) an unmarked pendant vertex $x$ and a second marked pendant vertex connected, via a marked edge, to a $\Komplete$-bag. In both cases, the vertex $x$ is called the \emph{special vertex} of $D$. 
\end{definition}

It follows that, if in the $\Komplete$-bag $K_t$ in $\dt_{G_2}$ has $t=3$ marked vertices, then Lines~\ref{line:C}--\ref{line:one-special} are executed. There, just one special vertex $w$ is removed from $X$ before returning this set as $\mu$-set of $G_2$. Note that, in case there are no special subgraphs in $\dd_{G_2}$, two vertices $w$ and $w'$ are removed from $X$ at Line~\ref{line:two-vertices}. 
The proof that this is the minimum number of vertices to remove from $X$ to get a mutually visible set, as well as the proof of Theorem~\ref{thm:main}, is omitted due to space.

Concerning the last example, $\dt_{G_3}$ is obtained from a graph $G_3$ in which its directed canonical decomposition contains just two opposite t-arrows $a$ and $b$. In the example, $\algo$ executes the block of Lines~\ref{line:opposite}--\ref{line:opposite-end}. This is because there are special subgraphs, and hence it is enough to remove just one special vertex from $X$ (Line~\ref{line:remove1} or~\ref{line:remove2}). In general, in case there are no special subgraphs, the algorithm eliminates two unmarked vertices from $X$ (one from $V_{G}^t(a)$ and the other from $V_{G}^t(b)$, cf. Lines~\ref{line:two-vertices2}--\ref{line:end-special2}).

\section{Correctness}\label{sec:correctness}

We provide here the correctness and complexity of $\algo$.

\begin{lemma}\label{lem:Va}
    Given a distance-hereditary graph $G$, if $a$ is a t-arrow, then $V_G^h(a)$ is the set of all unmarked vertices of $\dd^h_G(a)$.
\end{lemma}
\begin{proof}
Let $a=(u,v)$ and assume that there is a vertex $x\in V(G)$ in $\dd^h_G(a)$ not in $V_G^h(a)$, that is $x$ is not reachable from $v$ via an alternating path. Then the shortest path from $v$ to $x$ must contain at least a vertex that is the center of a $\Star$-bag. Consider the first of these vertices starting from $v$. If it is marked, then the corresponding marked edge is an arrow. Otherwise, it is a $\sigma$-vertex. In both cases, $a$ cannot be a t-arrow. 
\end{proof}

\begin{lemma}\label{lem:Vis}
    Let $G$ be a distance-hereditary graph. Let $x,y\in V(G)$ and let $P$ be the $x,y$-shortest path in $\dd_G$.  Let $\Sigma$ be the set of the internal $\sigma$-vertices of $P$ and let $A$ be the set of the arrows $a=(u,v)$ in $\dd_G$ such that only $u$ is in $P$. Vertices $x$ and $y$ are $X$-visible if and only if $\Sigma\cap X=\emptyset$ and for each $a\in A$ there exists a vertex $w\in V^h_G(a)$ such that $w\not \in X$.
\end{lemma}
\begin{proof}
($\Longleftarrow$) If $\Sigma$ and $A$ are both empty, $P$ is an alternating path and then $x$ and $y$ are adjacent in $G$, hence mutually visible. Otherwise, $P$ is not an alternating path and there are pairs of consecutive unmarked edges, that are both incident to a center of a $\Star$-bag. These centers are vertices in $\Sigma$, if they are $\sigma$-vertices, or tails of the arrows in $A$. Let $u_1, u_2,\ldots, u_k$ be the sequence of these vertices in $P$ from $x$ to $y$. Let us define $v_i$ as $u_i$ if $u_i$ is a $\sigma$-vertex or as a vertex $w\in V^h_G(a)\setminus X$ if $u_i$ is the tail of an arrow $a\in A$. Then $x=v_0, v_1, v_2,\ldots, v_k, v_{k+1}=y$ is a path in $G$ as there is an alternating path between $v_i$ and $v_{i+1}$, for each $i=0,1,\ldots, k$, in $D_G$ and then an edge between $v_i$ and $v_{i+1}$ in $G$ by Lemma~\ref{lem:alternating}. Call this path $Q$. Path $Q$ is an induced path in $G$ since there are no alternating paths connecting non-consecutive vertices $v_i$ in $D_G$ and then, since $G$ is a distance-hereditary graph, the path is a shortest path in $G$. So $Q$ is a shortest path in $G$ and no internal vertex of $Q$ is in $X$, then $x$ and $y$ are $X$-visible.

($\Longrightarrow$) Let us assume that $\Sigma\cap X\not =\emptyset$. Let $u$ be a vertex in $\Sigma\cap X$. Since $P$ pass through $u$ and by Lemma~\ref{lem:cut} $u$ is a cut-vertex in $G$, then each shortest path in $G$ between $x$ and $y$ pass through $u$. Then $x$ and $y$ are not $X$-visible. Assume now that there exists $a\in A$ such that  all the vertices in $V^h_G(a)$ are in  $X$. Consider any shortest path $Q$ between $x$ and $y$ in $G$. By Lemma~\ref{lem:alternating} each edge in $Q$ corresponds to an alternating path in $D_G$. Then sequences of alternating paths connect $x$ and $y$ in $\dd_G$ and since $P$ is unique in $\dd_G$ one of these alternating paths must pass throw $a$. Then a vertex of $V^h_G(a)$ is in $Q$. Hence $x$ and $y$ are not $X$-visible since $V^h_G(a)\subseteq X$.
\end{proof}

\begin{lemma}\label{lem:mu}
    Given a distance-hereditary graph $G$, Algorithm \A computes a mutual-visibility set for $G$.
\end{lemma}
\begin{proof}
Algorithm \A computes $\dd_G$ and a set $X$ as $V(G)$ minus the $\sigma$-vertices of $\dd_G$, that is minus the cut-vertices of $G$. By Lemma~\ref{lem:cut}, $X$ is a candidate to be a mutual visibility set. In the following, we analyse four cases depending on the number of t-arrows in $\dd_G$. 

\begin{itemize}
\item 
{\bf Case 1: $\dd_G$ has no t-arrows.} 
In this case Algorithm \A returns $X$. Let us show that $X$ is a mutual-visibility set. 
If $\dd_G$ has no arrows at all, by Lemma~\ref{lem:block}, $G$ is a block graph. A $\mu$-set for a block graph is given by all its vertices but its cut-vertices~\cite{DiStefano22}.  
Assume now that $\dd_G$ has arrows that are not t-arrows. Then, by definition of t-arrow, $V^h_G(a)$ contains a $\sigma$-vertex for each arrow in $\dd_G$. Then, since all the $\sigma$-vertices are not in $X$, each pair of vertices $x$ and $y$ are $X$-visible by Lemma~\ref{lem:Vis}.



\item 
{\bf Case 2: $\dd_G$ has one t-arrow or the t-arrows in $\dd_G$ are pairwise tail-connected.}  In this case Algorithm \A returns $X$ minus a vertex in $V_G^h(a)$ for each t-arrow $a$ in $\dd^h_G(a)$.
Lemma~\ref{lem:2arrows} ensures that for each arrow $a$ in $\dd_G$,  $V_G^h(a)$ has a vertex not in $X$, and then, by Lemma~\ref{lem:Vis}, any pair of vertices $x,y\in X$ are in mutual visibility.



\item 
{\bf Case 3: the t-arrows in $\dd_G$ are pairwise head-connected.}  By Lemma~\ref{lem:clique} all the heads of the t-arrows are all in a $\Komplete$-bag.
Algorithm \A returns $X$ minus an unmarked vertex $w$ of the $\Komplete$-bag, if it exists.
Then, since $w\in V_G^h(a)$ for each t-arrow $a$, by Lemma~\ref{lem:2arrows} also any arrow $b$ has an unmarked vertex in $V_G^h(b)$. By considering that the cut-vertices of $G$ are not part of $X$, by Lemma~\ref{lem:Vis} any pair of vertices $x,y\in X$ are in mutual visibility.

If all the vertices of the $\Komplete$-bag are marked Algorithm \A returns $X$ minus one or two vertices chosen among the vertices in $\dd_G^t(a)$ for some t-arrows $a$. Let $C$ be the set defined as in Line~\ref{line:C}. If $C$ does not contain a special subgraph, then Algorithm \A chooses two vertices from two subgraphs in $C$ and then for each arrow $a$ in $\dd_G$,  $V_G^h(a)$ has a vertex not in $X$. By Lemma~\ref{lem:2arrows} and Lemma~\ref{lem:Vis}, any pair of vertices $x,y\in X$ are in mutual visibility.  If $C$ contains a special subgraph, Algorithm \A removes from $X$ its special vertex and again for each arrow $a$ in $\dd_G$,  $V_G^h(a)$ has a vertex not in $X$.

\item 
{\bf Case 4: $\dd_G$ has exactly two opposite t-arrows.}  Call $a=(u,v)$ and $b=(v,u)$ the t-arrows. In case $\dd^h_G(a)$ and $\dd^h_G(b)$ are not special, Algorithm \A returns $X$ minus two vertices: a vertex $w_a$ in $V_G^h(a)$ and a vertex $w_b$ in $V_G^h(b)$.  
Note that $V^h_G(a)=V^t_G(b)$ and $V^h_G(b)=V^t_G(a)$ and by Lemma~\ref{lem:Va} these sets are the sets of all unmarked vertices of $\dd^h_G(a)$ and $\dd^h_G(b)$ respectively. 
Let $x,y$ be two vertices in $X$ and let $P$ be a shortest path connecting them in $D_G$.
 If $P$ pass on $u$ and $v$, then there exists an alternating path connecting $x$ and $y$ in $\dd_G$, hence they are adjacent in $G$ and hence in mutual visibility. If $P$ pass on $u$ only ($v$ only, resp.) $x$ and $y$ are in mutual visibility as shown in Case 2, whereas if $P$ pass neither on $u$ or $v$ $x$ and $y$ are in mutual visibility as shown in Case 1.

If at least one of $\dd^h_G(a)$ and $\dd^h_G(b)$ is special, say $\dd^h_G(a)$, then it is sufficient to remove only one vertex in $\dd^h_G(a)$ from $X$ to make all the vertices of $X$ in mutual visibility. In this case, Algorithm \A correctly returns $X$ minus a vertex $w_a$ in $V_G^h(a)$. 
\end{itemize}
\end{proof}

\begin{theorem}\label{thm:main}
     Given a distance-hereditary graph $G$, Algorithm \A computes a $\mu$-set for $G$ in linear time.
\end{theorem}
\begin{proof}
Lemma~\ref{lem:mu} shows that Algorithm \A computes a mutual visibility set $X$ for $G$.
Let us prove that $X$ is a $\mu$-set. By Lemma~\ref{lem:cut} we assume that $X$ does not contain any  vertex in the set $C$ of the cut vertices of $G$. If $\dd_G$ has no t-arrows, the set $X$ returned by algorithm $\A$ is given by all the vertices in $V(G)\setminus C$. Then $X$ is a $\mu$-set. If $\dd_G$ has only one t-arrow or tail-connected t-arrows,  $X$ is given by all the vertices in $V(G)\setminus C$ but one vertex for each $V_G^h(a)$, where $a$ is a t-arrow. Since $V_G^h(a)\cap V_G^h(b)=\emptyset$ for each pair $a,b$ of t-arrows, by Lemmata~\ref{lem:cut} and~\ref{lem:Vis} $X$ is a $\mu$-set. Indeed, if $V_G^h(a)\subseteq X$  for a certain t-arrow $a$ then there would be a pair of vertices in $D_G^t(a)$ not in mutual visibility. For the same reason, when $a$ and $b$ are opposite t-arrows, $X$ cannot contain one vertex in $\dd_G^h(a)$ and one vertex in $\dd_G^h(b)$ (unless $\dd_G^h(a)$ and $\dd_G^h(b)$ are special subgraphs of $\dd_G$). In this case $X$ is again a $\mu$-set. As discussed, if at least one of $\dd_G^h(a)$ and $\dd_G^h(b)$ is a special component, it is sufficient to remove only one vertex in $\dd_G^h(a)$ from $X$ in order to make it a $\mu$-set.

As for the last case where the t-arrows are head-connected, $X=V(G)\setminus (C\cup \{w\})$, where $w$ is an unmarked vertex in the $\Komplete$-bag containing all the heads of the t-arrows, if it exists. The removal of a vertex $w$ is needed otherwise pair of unmarked vertices in $\dd_G^t(a)$, for each t-arrow $a$ would not be in visibility. Then $X$ is a $\mu$-set. If such vertex does not exist, consider a t-arrow $a$. In order to make all the pair of vertices in $\dd_G^t(a)$ in mutual visibility, at least one vertex $w$ in $\dd_G^t(b)$, for some arrow $b\not = a$, must be not in $X$. Then all the pairs of vertices in $\dd_G^t(a')$, are in mutual visibility, for each $a'\not = b$. In order to make in mutual visibility all the pair of vertices in $\dd_G^t(b)$, when $\dd_G^t(b)$ is not special, we need to remove a further vertex $w'$ from $X$. Then $|X|=|V(G)\setminus C|-2$, and $X$ is a $\mu$-set. Finally, if $\dd_G^t(b)$ is special Algorithm \A chooses the special vertex $w$ of $\dd_G^t(b)$ and then no further vertex should be removed from $X$. Hence $|X|=|V(G)\setminus C|-1$ and $X$ is a $\mu$-set.

As for the computational complexity, $D_G$ can be computed in linear time (cf. Theorem~\ref{teo:dahlhaus_complexity})  as well as $\dd_G$ and $\dt_G$.
The set $S$ of all the $\sigma$ vertices of $\dd_G$ can be computed in linear time. The same time is needed for the set $A$ of all the t-arrows of $\dd_G$ that can be computed  by a visit in post-order of $\dt_G$ rooted to any vertex.


The cycle after Line~\ref{line:tails} can be computed in linear time since the computation of all the $V_G^h(a)$ for each $a\in A$ is bounded by the size of $\dd_G$.

If $A$ contains two opposite t-arrows, the algorithm checks if $\dd^h_G(a)$ is special, an operation that can be done in linear time.

The time to perform instructions at Lines~\ref{line:heads} and~\ref{line:K} is bounded by the size of the $\Komplete$-bag $K$, whereas the computation of $\dd^t_G(a)$ for each $a\in A$ at Line~\ref{line:C} is bounded by the size of the whole $\dd_G$, which is linear.
\end{proof}

%
\section{Conclusions}\label{sec:conclusion}
This work suggests some further research directions. The first natural open problem concerns, in general, the possibility of using the canonical decomposition also for computing the general position number $\gp(G)$ for a distance-hereditary graph $G$. To this end, perhaps the given concept of directed canonical decomposition, together with the given structural properties, could help.
If successful, one could try to use the canonical decomposition to compute both $\mu(G)$ and $\gp(G)$ for other classes of graphs for which the split decomposition has provided a characterization, such as circle graphs and parity graphs.






\begin{thebibliography}{10}
\expandafter\ifx\csname url\endcsname\relax
  \def\url#1{\texttt{#1}}\fi
\expandafter\ifx\csname urlprefix\endcsname\relax\def\urlprefix{URL }\fi
\expandafter\ifx\csname href\endcsname\relax
  \def\href#1#2{#2} \def\path#1{#1}\fi

\bibitem{DiStefano22}
G.~{Di Stefano}, Mutual visibility in graphs, Applied Mathematics and
  Computation 419 (2022) 126850.
\newblock \href {https://doi.org/10.1016/j.amc.2021.126850}
  {\path{doi:10.1016/j.amc.2021.126850}}.

\bibitem{CiceroneDK23}
S.~Cicerone, G.~{Di Stefano}, S.~Klavžar, On the mutual visibility in
  cartesian products and triangle-free graphs, Appl. Math. Comput. 438 (2023)
  127619.
\newblock \href {https://doi.org/10.1016/j.amc.2022.127619}
  {\path{doi:10.1016/j.amc.2022.127619}}.

\bibitem{CiceroneDKY23}
S.~Cicerone, G.~{Di Stefano}, S.~Klavžar, I.~G. Yero, Mutual-visibility in
  strong products of graphs via total mutual-visibility (2022).
\newblock \href {http://arxiv.org/abs/2210.07835} {\path{arXiv:2210.07835}}.

\bibitem{howorka:77}
E.~Howorka, A characterization of distance-hereditary graphs, Quarterly Journal
  of Mathematics 2~(28) (1977) 417--420.

\bibitem{bandelt/mulder:86}
H.-J. Bandelt, H.~M. Mulder, Distance-hereditary graphs, J. Comb. Theory Ser. B
  41~(2) (1986) 182--208.

\bibitem{peerj21}
S.~Cicerone, G.~{Di Stefano}, Getting new algorithmic results by extending
  distance-hereditary graphs via split composition, PeerJ Comput. Sci. 7 (2021)
  e627.
\newblock \href {https://doi.org/10.7717/peerj-cs.627}
  {\path{doi:10.7717/peerj-cs.627}}.

\bibitem{dam01}
S.~Cicerone, G.~{Di Stefano}, Graphs with bounded induced distance, Discret.
  Appl. Math. 108~(1-2) (2001) 3--21.

\bibitem{BrandstadtLS99}
A.~Brandst\"{a}dt, V.~B. Le, J.~P. Spinrad, Graph classes: a survey, SIAM,
  Philadelphia, PA, USA, 1999.

\bibitem{bouchet:88}
A.~Bouchet, Transforming trees by successive local complementations, Journal of
  Graph Theory 4 (1988) 195--207.

\bibitem{FPS19}
P.~Flocchini, G.~Prencipe, N.~{Santoro (Eds.)}, Distributed Computing by Mobile
  Entities, Current Research in Moving and Computing, Vol. 11340 of LNCS,
  Springer, 2019.
\newblock \href {https://doi.org/10.1007/978-3-030-11072-7}
  {\path{doi:10.1007/978-3-030-11072-7}}.

\bibitem{LunaFCPSV17}
G.~A.~D. Luna, P.~Flocchini, S.~G. Chaudhuri, F.~Poloni, N.~Santoro,
  G.~Viglietta, Mutual visibility by luminous robots without collisions, Inf.
  Comput. 254 (2017) 392--418.
\newblock \href {https://doi.org/10.1016/j.ic.2016.09.005}
  {\path{doi:10.1016/j.ic.2016.09.005}}.

\bibitem{PoudelAS21}
P.~Poudel, A.~Aljohani, G.~Sharma, Fault-tolerant complete visibility for
  asynchronous robots with lights under one-axis agreement, Theor. Comput. Sci.
  850 (2021) 116--134.
\newblock \href {https://doi.org/10.1016/j.tcs.2020.10.033}
  {\path{doi:10.1016/j.tcs.2020.10.033}}.

\bibitem{CiceroneFSN23}
S.~Cicerone, A.~{Di Fonso}, G.~{Di Stefano}, A.~Navarra, The geodesic mutual
  visibility problem for oblivious robots: the case of trees, in: 24th Int.
  Conference on Distributed Computing and Networking, {ICDCN} 2023, {ACM},
  2023, pp. 150--159.
\newblock \href {https://doi.org/10.1145/3571306.3571401}
  {\path{doi:10.1145/3571306.3571401}}.

\bibitem{ManuelK18}
P.~D. Manuel, S.~Klavzar, The graph theory general position problem on some
  interconnection networks, Fundam. Informaticae 163~(4) (2018) 339--350.

\bibitem{cunningham/edmonds:80}
W.~H. Cunningham, J.~Edmonds, A combinatorial decomposition theory, Canadian
  Journal of Mathematics 32 (1980) 734--765.

\bibitem{dahlhaus:jagm00}
E.~Dahlhaus, Parallel algorithms for hierarchical clustering and applications
  to split decomposition and parity graph recognition, J. Algorithms 36~(2)
  (2000) 205--240.
\newblock \href {https://doi.org/10.1006/jagm.2000.1090}
  {\path{doi:10.1006/jagm.2000.1090}}.

\bibitem{Courcelle06}
B.~Courcelle, The monadic second-order logic of graphs {XVI} : Canonical graph
  decompositions, Log. Methods Comput. Sci. 2~(2) (2006).
\newblock \href {https://doi.org/10.2168/LMCS-2(2:2)2006}
  {\path{doi:10.2168/LMCS-2(2:2)2006}}.

\end{thebibliography}


\end{document}